\documentclass[graybox]{svmult}

\usepackage{type1cm}        
\usepackage{makeidx}         
\usepackage{graphicx}        
\usepackage{multicol}        
\usepackage[bottom]{footmisc}

\usepackage{xcolor}
\usepackage[utf8]{inputenc}

\usepackage{bm}

\makeindex             

\begin{document}

\title*{Solving clustering as ill-posed problem: experiments with K-Means algorithm}
\author{Alberto Arturo Vergani}
\institute{Sant'Anna School of Advanced Studies \at The BioRobotics Institute, Pontedera (PI), Italy \at \email{albertoarturo.vergani@santannapisa.it}}
%
%
\maketitle


\abstract{In this contribution, the clustering procedure based on K-Means algorithm is studied as an inverse problem, which is a special case of the ill-posed problems. The attempts to improve the quality of the clustering inverse problem drive to reduce the input data via Principal Component Analysis (PCA). 
Since there exists a theorem by Ding and He that links the cardinality of the optimal clusters found with K-Means and the cardinality of the selected informative PCA components, the computational experiments tested the theorem between two quantitative features selection methods: Kaiser criteria (based on imperative decision) \textit{versus} Wishart criteria (based on random matrix theory). The results suggested that PCA reduction with features selection by Wishart criteria leads to a low matrix condition number and satisfies the relation between clusters and components predicts by the theorem. The data used for the computations are from a neuroscientific repository: it regards healthy and young subjects that performed a task-oriented functional Magnetic Resonance Imaging (fMRI) paradigm.}

\section{Introduction}


\subparagraph{The clustering procedure}
Clustering, or partitioning, is a computational methodology used in many fields to draw knowledge on the distribution of raw data and the informative patterns that are within \cite{Halkidi2001, Li2016}. The main objective of clustering is to collect into groups, i.e., clusters, a given input data taking account a similarity measures predefined, e.g., algebraic distance \cite{Pakhira2004}. Formally, let $ X = \{ x_{1}, x_{2}, \dots, x_{N} \} $ be a given input dataset with $ N $ points, $ K $ is the number of clusters, also known as patterns, within the data. The aim of clustering is to compute a partition matrix $ U(X) $ of the data in order to define a partition $ C = \{ C_{1}, C_{2}, \dots, C_{k} \} $, with the criteria that similar points have high similarity and dissimilar points have high dissimilarity, i.e. the similar points are close as possible, while the dissimilar points are far as possible. The partition matrix is denoted as $ U = [\mu_{ij}] $, $ 1 \leq i \leq K $, $ 1 \leq j \leq N $, where $ \mu_{ij} $ is the grade of membership of point $ x_{j} $ to cluster $ C_{i}(i=1, \dots, K) $.

It is important to note that clustering can be performed in two radically different ways: crisp and fuzzy.
\begin{itemize}
	\item The crisp clustering uses a logical framework based on classic sets theory, i.e. any one point of a given dataset belongs to only one class of clusters. It means that there exists a pattern, and only one pattern, for each points. Formally, the membership value $ \mu_{ij} $ has value equal to 1 if $ x_{j} \in C_{i} $, otherwise  $ \mu_{ij} = 0 $.
	\item The fuzzy clustering is based on the fuzzy sets theory, i.e., a point may belong to one or more clusters with a certain grade of membership. The partition matrix $ U(X) $ is represented as $ U = [\mu_{ij}] $, where $ \mu_{ij} \in [0, 1] $. Therefore, the crisp clustering is a special case of fuzzy clustering if the membership of a point to a specific cluster is either 0 or 1. Furthermore, there exist methods to hard the fuzzy clustering, e.g., selecting for each points the cluster with associated the largest grade of membership. 
\end{itemize}

\subparagraph{The K-Means clustering algorithm}
One of the most important crisp partition method is the K-Means algorithm. It was introduced by MacQueen in the Sixties \cite{MacQueen1967} in the context of signal processing as method for vector quantization. It works with $ k $ predefined clusters used to partition $ n $ points, i.e., the observations, with the inclusion criteria based on their similarity with a prototype. In features space, the position of the $ k $ centroids is used to evaluate the membership of each points, i.e., one point is within a cluster in which has nearest mean. 

More formally, given a set of observation $ (\bm{X}_{1}, \bm{X}_{2}, \dots, \bm{X}_{n} ) $, where each observation is point in a $ d $ -dimensional features space, the k-means clustering partitions the $ n $ observations into $ k $ classes $ C=\{C_{1}, C_{2}, \dots, C_{k}\} $ according to the minimization of the variance within the cluster using an iterative procedure that minimizes the squared error of the following objective function:

\begin{equation} \label{eq:kmeans}
J = \sum_{i=1}^{n} \sum_{c=1}^{C} ||X_{i(c)} - A_{c}||^{2}
\end{equation}

where $ ||...||^{2} $ is a specified distance function, $ n $ are the data points, $ k $ are the clusters, $ X_{i(c)} $ are the points in cluster $ c $, $ A_{c} $ are the centroids of clusters $ c $.

After that clustering algorithm computed the classes, it is necessary to evaluate if the procedure found an optimal solution. In order to do this, there are several measures able to state if clustering has been done properly. One of them is the separation measure proposed by Davies and Bouldin for crisp clustering (DB index) \cite{Davies1979}, which is defined as the system-wide average similarities of each clusters. Using the same notation of the authors, the index is formally:

\begin{equation} \label{eq:sdb1}
\overline{R}=\frac{1}{C}\sum_{i=1}^{C} R_{i} 
\end{equation}

where $ C $ is the number of clusters and $ R_{i}  $ is the $ max_{j \neq i}\{ R_{i,j} \} $. $ \overline{R} $ is the system-wide average of the similarity measures of each cluster with its most similar cluster; therefore, in relation to the algebraical properties that shape the DB index, the better the clusters, the lesser the average similarity. With  $ R_{i,j} $ is indicated the \textit{ratio} between the average similarities of each clusters with its most similar cluster and the distance between the centroids of the clusters. More in details: given the average similarities $ S_{i} $ and $ S_{}j $ of the cluster $ i $ and the cluster $ j $, and given $ M_{i,j} $ the distance between the centroids of the cluster $ i $ and $ j $, $ R_{i,j} $ is defined as:

\begin{equation} \label{eq:sdb2}
R_{i,j} = \frac{S_{i} + S_{j}}{M_{i,j}}
\end{equation} 

where, formally, $ S_{i} $ is :

\begin{equation} \label{eq:sdb3}
S_{i} = \Bigg\{ \frac{1}{T_i} \sum_{j=1}^{T_{i}} ||X_{j} - A_{i} ||^{q} \Bigg\}^{\frac{1}{q}}
\end{equation} 

in which $ T_{i} $ is the number of vectors in the clusters $ i $, $ X_{j} $ are the vectors in $ i $ and $ A_{i} $ is the centroid of the cluster $ i $. When $ q=1 $, $ S_{i} $ is the average Euclidean distance of vectors to the centroids of $ i $, whereas if $ q=2 $, $ S_{i} $ is the standard deviation of metric about the samples in cluster in relation with its centroid. $ S_{j} $ is defined equivalently. Instead, the distance between centroids is formalized as: 

\begin{equation} \label{eq:sdb4}
M_{i,j} = \Bigg\{ \sum_{k=1}^{N} ||a_{ki} - a_{kj} ||^{p} \Bigg\}^{\frac{1}{p}}
\end{equation}

where $ a_{ki} $ is the \textit{k}th component of the \textit{n}-dimensional vector $ a_{i} $, which is the centroid of cluster $ i $. $ M_{i,j} $ is the Minkowski metric of the centroids $ i $ and $ j $; than, if $ p=1 $ is the city block/Manhattan distance, whereas if $ p=2 $ is the Euclidean distance between centroids; note that if $ p=q=2 $, $ R_{i,j} $ is the Fisher similarity measure computed between cluster $ i $ and cluster $ j $.

\subparagraph{Clustering as an ill-posed problem}

According to Hadamard's definition \cite{hadamard1902problemes} a problem is said to be ill-posed if at least one of the following criteria is not satisfied: 1) a solution exists, 2) the solution is unique and 3) its behaviour depends continuously with the data. \textit{Viceversa}, if they are completely satisfied, the problem is defined as well-posed.

A subclass of ill-posed problem are the inverse problems, i.e., an approach that wants to discover causes of an observed effect, e.g., in helioseismology that one wants to determine the structure of the sun by measurements from earth or space, in image restoration that one wants to determine unavailable exact image from an available contaminated version, or in medical imaging where one wants to discover physiological causation related to experimental paradigm. 

Clustering is an ill-posed problem since it is lacking at least of the property to have a unique solution, because there are many possible clusters that can make interesting classifications of the input data. Specially, clustering is an inverse problem as its aim is to find patterns, i.e., causes/categories that explain the experimental data measures.

According to Bezdek \cite{Bezdek1981, Bezdek2016}, there are other issues concerning the clustering analysis: 1) the \textit{a priori} assumption that datasets have clusters, 2) the computation method to find clusters, 3) the evaluation procedure to verify the clusters found.

Therefore, merging the requirements to solve the illness defined by Hadamard and the issues proposed by Bezdek, the inverse clustering problem could be address with the following strategy:

\begin{itemize}
	\item \textit{the existence of a solution must be an assumption}, i.e., if the goal is to cluster a dataset, the experimentalist need to assume that there are at least two clusters;
	\item \textit{the uniqueness of a solution must be a decision}, according to some optimal evaluation methodology, e.g., using a measure able to define what is the optimal clustering result, e.g., using the Davies-Bouldin index before mentioned;
	\item \textit{the matrix condition number must be lowered}. Usually the condition number of experimental matrix is very high. Therefore, it is necessary to improve the quality of the matrix to have a more stable solution, i.e., small changes in the matrix should be associated to small changes in the results. The matrix condition number decreases using signal processing (e.g., Empirical Mode Decomposition) and data reduction techniques (e.g., Principal Component Analysis).
\end{itemize}

\subparagraph{The relation between K-Means and PCA}
The principal component analysis (PCA)\footnote{To look more details about the mathematical fundation of PCA and its several applications see the book by Jollife \cite{Jolliffe2011}.} was ideate by Pearson in the first decade of XX century \cite{Pearson1901} and independently discovered by Hotelling in the Thirties \cite{Hotelling1933}.
It is an unsupervised linear transformation technique that is widely used in many fields, e.g., science, engineering, biomedical domains, with the specific aim to reduce the dimensionality of a system \cite{raschka2015python}.

Basically, PCA identifies patterns in data based on the correlation matrix of the features. The analysis of the principal component aims to find the direction of the maximum variance in high dimensional data and projects it onto a new subspace with equal or fewer dimensions than the original system. Therefore, the principal components are the orthogonal axis of the new subspace and they can be understood as the directions of the maximum variance given the constraints that the new features' axes are orthogonal to each other \cite{raschka2015python}.

More formally, given $ Y = (\bm{y}_{1}, ..., \bm{y}_{n}) $ as the standardize data matrix, the covariance of $ Y $, ignoring the factor $ 1/n $, is $ \Sigma_{i} = (\bm{x}_{i} - \overline{\bm{x}}) (\bm{x}_{i} - \overline{\bm{x}})^{T} = YY^{T} $. The principal directions $ \bm{u}_{k} $ and the principal components $ \bm{v}_{k} $ are eigenvectors satisfying
\begin{equation}
YY^{T}\bm{u}_{k} = \lambda_{k} \bm{u}_{k} 
\end{equation}

\begin{equation}
Y^{T} Y \bm{v}_{k} = \lambda_{k} \bm{v_{k}} 
\end{equation}

\begin{equation}
\bm{v}_{k} = Y^{T} \bm{u}_{k} / \lambda_{k}^{1/2}
\end{equation}

Elements of $ \bm{v}_{k} $ are the projected values of data points on the principal directions $ \bm{u}_{k} $.

According to the work by Ding and He there is relation between K-Means clustering algorithm and the Principal Component Analysis (PCA). In their seminal paper \cite{ding2004k}, they proved that principal components are actually a continuous solution of the cluster membership indicators for the K-Means clustering method. It means that the PCA dimension reduction automatically performs data clustering according to the K-Means objective function $ J(K) $ with $ K $ clusters. Also, they linked the number of clusters $ K $ with number of principal directions found by PCA, i.e., the clustering structure is embedded in the  $ K - 1 $ components. More formally, they proposed the following theorem \cite{ding2004k}\footnote{The proof of this theorem is the direct application of the well-known theorem by Ky Fan to optimization problem \cite{fan1949theorem}}

\begin{quotation}\label{quot:ding}
	When optimizing the K-Means objective function, the continuous solution for the transformed discrete cluster membership indicator vectors $ Q_{K-1} $ are the $ K-1 $ principal components $ Q_{K-1} = (\bm{v}_{1}, ..., \bm{v}_{K-1} ) $. The functional\footnote{$ J{K} $ is K-means functional to be optimize minimizing the sum of squared error $ J{K} $ (see Equation \ref{eq:kmeans}} $ J{K} $ satisfied the upper and lower bounds
	
	\begin{equation}
	n \overline{\bm{y}^{2}} - \sum_{k=1}^{K-1} \lambda_{K} <  J{K} < n \overline{\bm{y}^{2}}
	\end{equation}
	
	where $ n \overline{\bm{y}^{2}}  $ is the total variance and $ \lambda_{k} $ are the principal eigenvalues of the covariance matrix $ YY^{T} $
	
\end{quotation}
Practically, if the correct number of clases is three in a selected dataset (e.g., the Fisher's Iris dataset \cite{Fisher1936}), the most discriminative subspace is shaped by the first two eigenvector, because this PCA subspace is particularly effective for the K-Means clustering.  

\subparagraph{Objectives}
Since the above framework linked the clustering to an inverse problem, the main objective of this work is to frame the clustering procedure within the ill posed problem conceptualization, in order to solve some issues that are related with it. 

In consideration that the existence of a solution is an assumption and that the uniqueness of the solution is a decision (taken with clustering evaluation index), the continuity property that data needs to have with the solution may be checked and - nevertheless - improved.

Therefore, the first step is compute the matrix condition number that may indicates how much is ill-conditioned the data matrix. The second step is to test the theorem of Ding and He in order to find experimentally the relation between the cardinality of the optimal clusters and the cardinality of the informative components. 

The datasets used to do this operations derive from experiments coming from functional Magnetic Resonance Imaging (fMRI), which is a noninvasive technique for investigate the activation of the human brain. Precisely, 
\begin{itemize}
	\item during a functional experiments, a series of brain images are acquired and the variations in the measured signals between images are used to infer the brain functionality related to specific task \cite{Lindquist2008}; 
	\item the fMRI data consists of a number of uniformly spaced volume elements, i.e. the voxels, which is a partition of the brain in equally size boxes and the intensity of each voxel represents for that area the spatial distribution of the nuclear spin density \cite{Lindquist2008};
	\item the changes in the brain hemodynamic, as a reaction to neural activity, impact on the MR signal, therefore the changes in voxel intensity across time are an indirect measure of the cells activations and this measurement can be used to infer the spatial and the temporal information regard some brain functionality \cite{Ogawa1992};
	\item the analysis of the fMRI data regard multiple objectives: the localization of the regions involved in the task, the determination of the distributed networks linked to brain functions, the predictions about physiological and pathological brain states \cite{Lindquist2008}.
\end{itemize}

\section{Experimental Procedures}

\subsection{The fMRI datasets}
The dataset selected is the one proposed by Wakeman and Henson \cite{Wakeman2015} from \textbf{openfmri} repository (https://openfmri.org/dataset/ds000117). This dataset regards a study where 16 healthy and young subjects did a task-based paradigm, which consists of a sequence of task-blocks (when subjects perform cognitive tasks) alternated with rest-blocks (when subjects did anything). In particular the experiment is a fMRI Face Recognition Task, which is a repeated sequences of familiar (famous) faces, unfamiliar (non-famous) faces and scrambled faces. Each stimulus type is supposed to elicit a different physiological brain response. 

Wakeman and Henson used the Statistical Parametric Maps (SPMs) \cite{friston1994statistical} to localize the brain activity associated to each type of stimulus, which is a model-based approach to analyse fMRI data. In this study, instead, the goal is to investigate with clustering the functional patterns in order to find time-related brain activities of interest, i.e., to detect collective regional activations associated to specific time points during the experimental blocks, which are also repeated along the experimental paradigm. Since temporal-related brain patterns are informative about the similarity between physiological statuses, which is a neuronal common response to specific task-block, the objective of clustering is this study is to find, specifically, fMRI temporal patterns. 

\subsection{Computational Steps}	
The computational procedures regard the steps of 1) data elaboration and reduction, 2) condition matrix evaluation and 3) the test of Theorem by Ding and He.

\begin{enumerate}
\item \textbf{Data elaboration and reduction:} 

\begin{itemize}
	\item \textit{Dataset RAW}. The original datasets processed with standard image processing tools \cite{Jenkinson2012}: such as spatial and temporal filtering, motion correction, standard registration (with MNI 152 reference), time-series extraction according to the meaning of Harvard-Oxford Atlas with 96 lateralized labels \cite{Frazier2005, Makris2006}. The results is a matrix $ X \in \Re^{nxm} $ where $ n = (1, ..., 96) $, i.e., 96 ROIs, and $ m = (1, ..., 208) $, i.e., 208 time points (aka: brain volumes or experimental blocks). Therefore, each dataset is a matrix contain 96 time series vectors with 208 samples of the brain signal.
	
	\item \textit{Datasets EMD}. The dataset EMD are the datasets RAW with also the computation of the Empirical Mode Decomposition \cite{wang2010intrinsic,huang1998empirical}, which it decomposes the ROI's signal into a $ k $ number of Intrinsic Mode Function (IMF) and a residual process $ r_{k}(t) $, taking finally the 1st $ k $ IMF for each ROIs.
	
	\item \textit{Dataset PCA-K}. The datasets PCA-K are the datasets EMD reduced via Principal Component Analysis selecting feature components according Kaiser criteria \cite{kaiser1960application}. The Kaiser rule allows to define informative the eigenvectors associated with eigenvalues greater than one, i.e., $ \lambda \geq 1 $. 
	
	\item\textit{Datasets PCA-W.}. The datasets PCA-W are the datasets EMD reduced via Principal Component Analysis selecting features components according the Wishart criteria, which is based on the Wishart distribution\footnote{It is noted also as Marchenko-Pastur distribution \cite{marvcenko1967distribution}.}. This rule is part of the Random Matrix Theory \cite{mehta2004random} framework. It allows to describe the eigenvalue density of the empirical correlation matrix of independent identically distributed (iid) Gaussian variables. Wishart distribution is given by 
	
	\begin{equation}
	\rho (\lambda) = \frac{1}{2 \pi r \lambda} \sqrt{(\lambda_{+} - \lambda)(\lambda - \lambda_{-})}
	\end{equation}
	
	The parameter $ r = N / T$ is the rectangular \textit{ratio} between the number of variables $ N $ and the sample size in the data matrix $ T $. According to the above equation, eigenvalues of the correlation matrix are defined in a finite support $ \lambda \in [\lambda_{-},\lambda_{+}] $ and the end-points are $ \lambda_{\pm} = (1 \pm \sqrt{r})^{2} $\footnote{Note that if $ r \longrightarrow 0 $, than $ T $ is much larger than $ N $ and the $ \rho(\lambda) $ became a $ delta $ function. When $ r = N/T $ is finite, the eigenvalue density is smeared. However, when $ N=T $, $ r=1 $ is critical because the $ \lambda_{-} = 0 $.}. With $ N = 96 $ and $ T = 208 $, $ \lambda_{\pm} = \pm 2.8 $. 
	Therefore, the usage of the Wishart criteria means to include eigenvectors that have eigenvalues greater than $ \lambda_{+} = 2.8 $.  The components associated with eigenvalues that are less than $ \lambda_{+} $ have to be define as non-informative or noisy components, since they are within the support of the eigenvalues associated to random correlations.

\end{itemize}

\item \textbf{Evaluation the matrix condition number}. A mathematical problem with low condition number is said well-conditioned (i.e., $ K(A) = 1 $), while if it has high condition number is said bad-conditioned (i.e., $ K(A) >> 1 $; e.g., singular matrices have $ K(A) = \infty $). Specifically, the definition of condition number used for this operation is the 2-norm condition number for inversion of the matrix $ A $, i.e., 

\begin{equation}
	K(A) = \frac{\lambda_{max}}{\lambda_{min}}
\end{equation} 
where $ \lambda_{max} $ is the largest eigenvalue of the matrix $ A $ and $ \lambda_{min} $ is the smallest eigenvalue of the matrix $ A $ \cite{quarteroni2014scientific}

\item \textbf{Test the Theorem of Ding and He}, i.e., the existence of a relation between the optimal clustering and the informative components \cite{ding2004k}. The exact rule $ K - 1 = C $, with $ K $ the correct number of clusters and $ C $ the number of components that allow to find them, it has been relaxed with the following heuristic rule\footnote{The heuristic rule is motivated by the explorative value that have clustering studies. Often, the patterns investigated are complex structures that need the usage of several methodological comparisons and human expertise to correctly understand results.}:
\begin{equation} \label{eq:eur}
K \sim C
\end{equation}
In the context of the clustering inverse problem, the $ K $ clusters are not given and they have to be computed and selected according to some clustering methods and evaluation metrics. In this specific case, the clustering method used is the K-Means algorithm and the evaluation procedure the Davies-Bouldin Index. Therefore, the correct $ K $ (that is unknown in the inverse clustering problem case) is substituted with the estimated $ K^{*} $ and the equation \ref{eq:eur} becomes
\begin{equation} \label{eq:ast}
K^{*} \sim C
\end{equation}
If the similarity relation between $ K $ and and $ C $ holds, also the similarity between $ K^{*} $ and $ C $ need to be consistent. Therefore, once the optimal clustering partitioning is computed, also the correct number of components are estimated. The components $ C $ are the subspace of the total eigenvectors, after PCA reduction, that permit to have an effective clustering partition. The selecting criteria regarding the components that has the most likeness with the estimated $ K^{*} $ is the selecting criteria able to find the correct number of clusters. In PCA reduction, the selection of the relevant components has done with two main quantitative criteria: the Kaiser method and the Wishart method. 

Given that the similarity between optimal clusters number and the informative components number holds, the criteria that estimates the components number more similar to the optimal clusters number is the criteria that allows the K-Means to cluster better the dataset. The comparison between optimal clustering and informative components attempts to validate which of the following statements is true:
\begin{equation} \label{eq:kerualk}
K^{*} \sim C^{K}
\end{equation} 
or,
\begin{equation} \label{eq:kequalcw}
K^{*} \sim C^{W}
\end{equation} 

where $ C^{K} $ are the components selected with the Kaiser criteria and $ C^{W} $ the components selected with the Wishart criteria. 
Since only one of the above proposition could be true, the one that is true is associated with the criteria to select components able to find the optimal clustering partition.

In results section will be presented the following outcomes: a) the relation between the condition number of the matrices obtained with the different preprocessing (RAW, EMD, PCA-K, PCA-W); b) the test of Theorem by Ding and He using the relaxed version of the rule that links the correct clusters number to the informative components (validating the equations \ref{eq:kerualk} and \ref{eq:kequalcw}). 

\end{enumerate}

\section{Commented Results}

\subparagraph{The variation of condition number between processed and reduced datasets} The Tab \ref{tab:condition_number} shows the condition number computed for the processed and reduced datasets. RAW datasets are the original fMRI datasets preprocessed with classical signal processing. EMD datasets are the RAW datasets with also the computation of the Empirical Mode Decomposition (taking the 1st mode). The PCA-K is the EMD datasets reduced with PCA having components selected by  Kaiser criteria. The PCA-W is the EMD datasets reduced with PCA having components selected by Wishart criteria. The table communicates a specific result: the lower condition number is associated to the dataset reduced via PCA with features selected according to Wishart criteria. Furthermore, this results indicates that the third requisite to be a well posed problem (the continuity properties of the input in relation to the output) that ill problems do not have, could be addressed reducing the original dataset (that have the higher condition number) and selecting the components according to a Wishart criteria that discriminates the informative components, excluding the ones associated to random correlations and noise.

\begin{table}[tbph!] 
	\centering
	\begin{tabular}{c|c|c|c|c}
		SUBJECTS & RAW & EMD & PCA-K & PCA-W \\ \hline \hline
		1 & 34431.68 & 77.65 & 6.69 & \textbf{3.83} \\ \hline
		2 & 43075.31 & 80.84 & 7.21 & \textbf{3.85} \\ \hline
		3 & 36633.43 & 63.31 & 5.85 & \textbf{3.34} \\ \hline
		4 & 37317.32 & 73.17 & 6.52 & \textbf{4.02} \\ \hline
		5 & 47175.60 & 45.92 & 5.50 & \textbf{2.95} \\ \hline
		6 & 36776.46 & 34.42 & 5.03 & \textbf{2.95} \\ \hline
		7 & 41331.35 & 43.07 & 5.15 & \textbf{2.86} \\ \hline
		8 & 39576.56 & 38.25 & 4.80 & \textbf{2.73} \\ \hline
		9 & 39511.79 & 60.98 & 4.95 & \textbf{3.04} \\ \hline
		10 & 40831.24 & 209.84 & 6.99 & \textbf{4.00} \\ \hline
		11 & 45199.41 & 50.24 & 6.00 & \textbf{3.51} \\ \hline
		12 & 35033.68 & 46.98 & 5.60 & \textbf{3.17} \\ \hline
		13 & 44149.23 & 35.99 & 4.82 & \textbf{2.79} \\ \hline
		14 & 43850.81 & 35.15 & 5.22 & \textbf{3.10} \\ \hline
		15 & 47698.72 & 33.38 & 4.40 & \textbf{2.52} \\ \hline
		16 & 34316.43 & 79.88 & 6.29 & \textbf{3.67} \\ 
	\end{tabular}
	\caption{Condition number computed for the processed and reduced datasets. RAW datasets are the original fMRI datasets preprocessed with classical signal processing. EMD datasets are the RAW datasets with also the computation of the Empirical Mode Decomposition (taking the 1st mode). The PCA-K is the EMD datasets reduced with PCA having components  selected by  Kaiser criteria. The PCA-W is the EMD datasets reduced with PCA having components selected by Wishart criteria} 
	\label{tab:condition_number}
\end{table}

\begin{figure}[tbph!]
	\centering
	\includegraphics[width=1.2\linewidth]{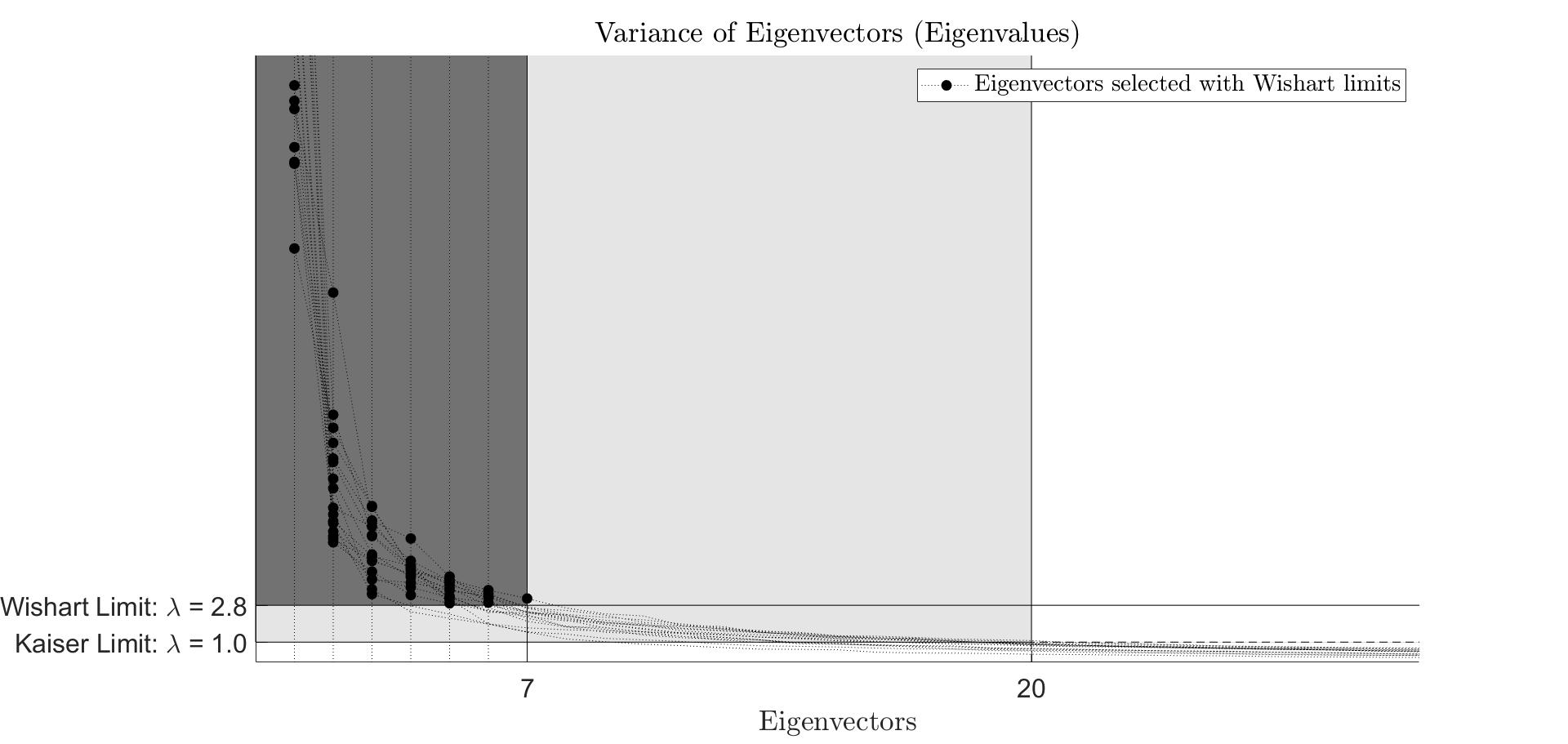}
	\caption{The figure shows the ordered variance of eigenvectors, i.e., the eigenvalues, about the all subjects. The eigenvalues greater than Kaiser criteria ($ \lambda=1 $) have associated the first twenty eigenvectors. Instead, the eigenvalues greater the Wishart limit ($ \lambda= 2.8 $) have associated the only the first seven eigenvectors.}
	\label{fig:variance-eigenvectors}
\end{figure}

\subparagraph{The variance of eigenvectors and their selection via Kaiser and Wishart methods} The Fig \ref{fig:variance-eigenvectors} shows the variance of eigenvectors, i.e., the eigenvalues, about all subjects. The eigenvalues greater than the Kaiser criteria ($ \lambda=1 $) have associated the first twenty eigenvectors. Instead, the eigenvalues greater the Wishart limit ($ \lambda_{max} = 2.8 $) have associated only the first seven eigenvectors. These results indicated different cardinality about the vectorial subspace used to compute the clusters with K-Means algorithm. The Kaiser criteria is more comprehensive than the Wishart criteria: it means that PCA data reduction followed by the features selection according with them gives different reduced datasets: the datasets reduced according to Kaiser criteria have more dimensions than the ones reduced according to Wishart criteria. 

\begin{figure}[tbph!]
	\centering
	\includegraphics[width=1.2\linewidth]{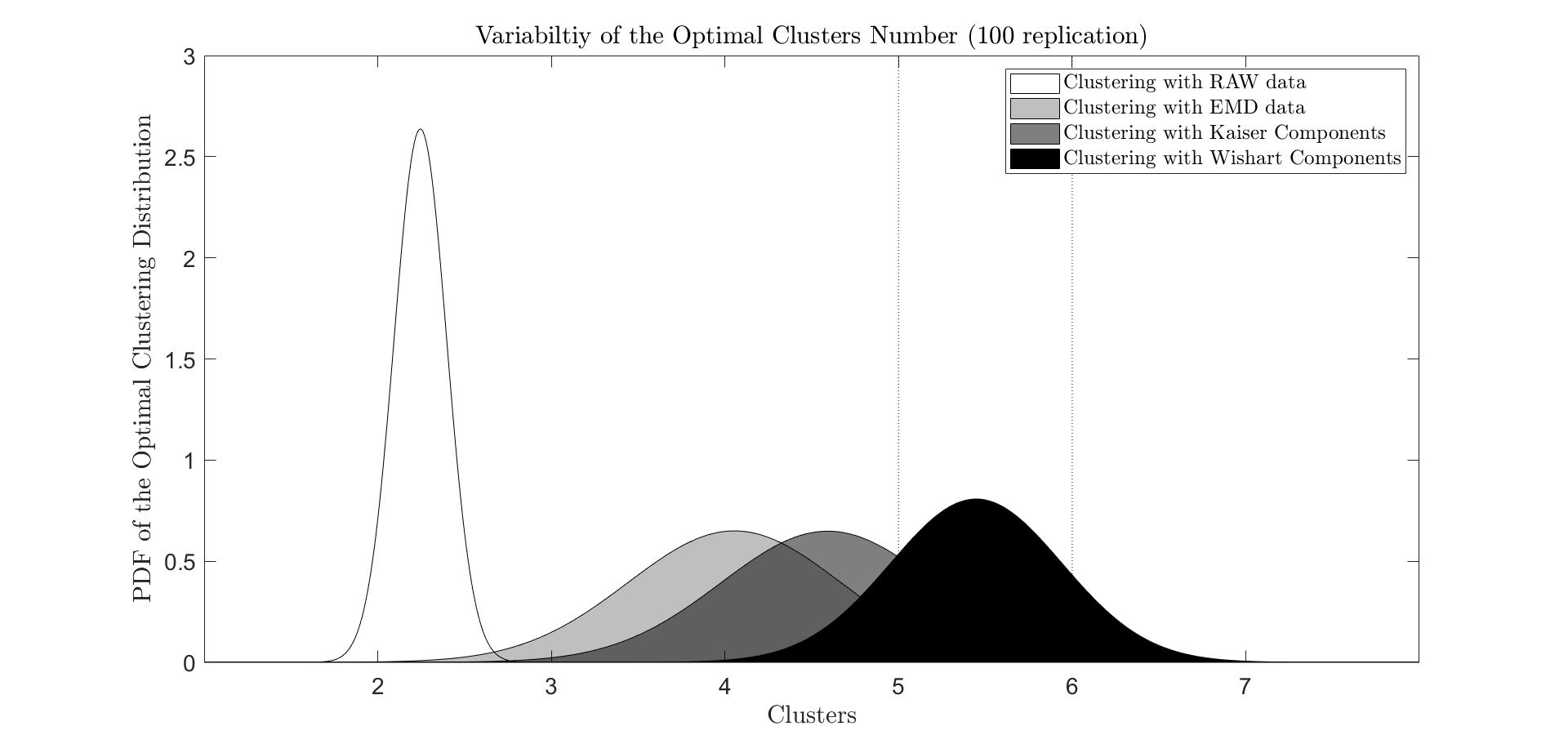}
	\caption{The figure represents the optimal cluster number for all the subjects using RAW data (white-empty shape), EMD data (soft gray shape), data reduced with PCA and components selection according to Kaiser criteria (hard gray shape) and data reduced with PCA and components selection according to Wishart criteria (black shape). The number of clusters found with dataset reduced with selected components according to Wishart criteria are very similar the number of the components used (see Fig \ref{fig:components-variations}).}
	\label{fig:cluster-variation}
\end{figure}

\begin{figure}[tbph!]
	\centering
	\includegraphics[width=1.2\linewidth]{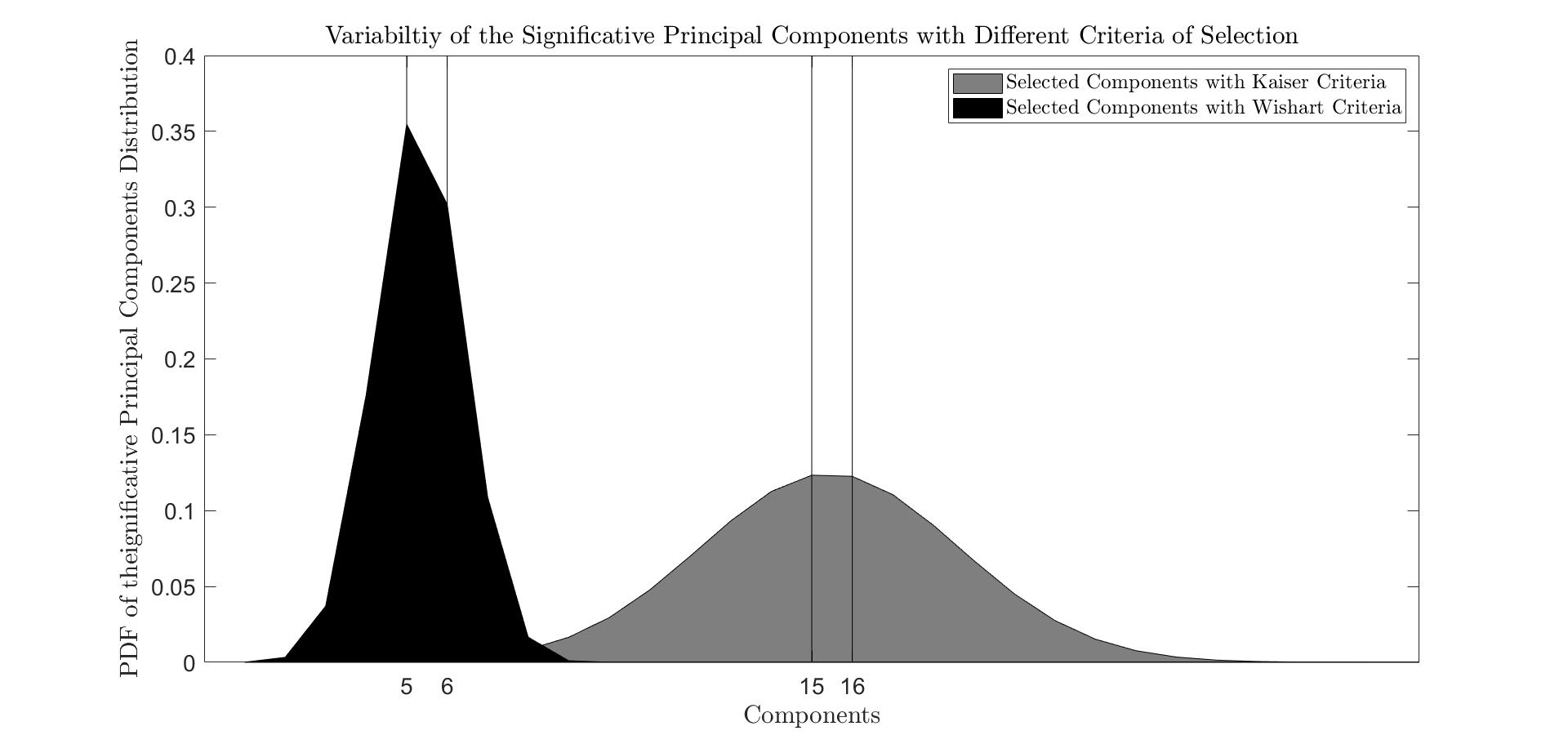}
	\caption{The figure shows the number of components selected with Kaiser criteria (hard gray shape) and the ones selected with Wishart criteria (black shape). The number of the components selected with Wishart criteria are very similar the number of clusters found used them as feature vectors (see Fig \ref{fig:cluster-variation}) }
	\label{fig:components-variations}
\end{figure}

\subparagraph{The variation of the optimal clustering between processed and reduced datasets} The Fig \ref{fig:cluster-variation} represents the optimal cluster number for all the subjects using RAW data (white-empty shape), EMD data (soft gray shape), data reduced with PCA and components selection according to Kaiser criteria (hard gray shape) and data reduced with PCA and components selection according to Wishart criteria (black shape). The different datasets lead to some heterogeneity into clustering results, i.e., the classes found with RAW datasets are few (2/3 patterns), whereas the classes found with datasets reduced according the Wishart components are many (4/7 patterns). The classes found with EMD datasets and with the datasets reduced according to Kaiser criteria are in the middle range of the white and black extremes. It is interesting to note that the number of clusters found with dataset reduced by selected components according to Wishart criteria are very similar to the number of the components themselves (see Fig \ref{fig:components-variations}).

\subparagraph{The variation of selected components via Kaiser and Wishart methods}
The Fig  \ref{fig:components-variations} shows the number of components selected with Kaiser criteria (hard gray shape) and the ones selected with Wishart criteria (black shape). There is an evident variation between the informative components selected with the two criteria with all subjects. Globally the components selected with Wishart criteria are less than the components selected with Kaiser criteria. 

Comparing the results of Fig\ref{fig:components-variations} with the results of Fig \ref{fig:cluster-variation}, it emerges the relation between the clusters number and the components number: the components selected with Wishart criteria are very similar with the number of clusters found using them as input features, whereas the components selected with Kaiser criteria are very different with the number of clusters found using them as input features.

Therefore, from these outcomes, it is possible to decide which of the relation \ref{eq:kerualk} and  \ref{eq:kequalcw} is true: considering the high similarity between the components selected according the Whisart criteria and the patterns found with K-Means clustering, it is reasonable to state that, from an experimental point-of-view, there is an heuristic link between the cardinality of the optimal clusters computed with K-means and with the cardinality of the informative eigenvectors subspace used as reduced input features, i.e, the relation \ref{eq:kequalcw} is experimentally true.

%

%

\section{Conclusions}
Since the clustering is an ill posed problem, some strategies have to be taken: 1) the existence of a solution need to be an assumption; 2) the uniqueness of solution may be a decision according a clustering evaluation criteria, which is able to measure the optimality of the clustering outcomes; 3) the continuity relation from input to output has to be manage with lowering the matrix condition number. The assumption about the existence of a solution is a classical assumption in the clustering theory domain (see Bezdek et al \cite{Bezdek2016}). The criteria to select only the optimal clustering solution is the Davies-Bouldin separation measure. The way to manage the condition number of the data matrix is the reducing procedure via Principal Component Analysis, and selecting the feature components according the Kaiser criteria and the Wishart criteria.

Furthermore, clustering and component analysis are not disjoint algorithms, but they shared hidden mathematical structures explored by Ding and He \cite{ding2004k}, which proposed a theorem that links the number of clusters with the number of components useful to get the optimal clustering: if one knows the true clusters, than also one knows the components space able to reduce the original dataset. However, in this computational study, the goal was to find the optimal number of clusters, managing the issues related to ill posed problem. Therefore, there is no knowledge  \textit{a priori}  about the true clusters number.

The theorem proposed by Ding and He \cite{ding2004k} predicts components from true clusters, than the same theorem could predict components from optimal clusters, which are clusters estimated as possibly true. According to this extension, the predicted components from the estimated optimal clusters need to be the most informative collection of eigenvectors able to reduce the original datasets. 

From the above propositions, the main computational results of this study could be formulate as the following statements:

\begin{quotation}
1) Data reduction via Principal Component Analysis with features selection according to Wishart criteria lowered the matrix condition number.
\end{quotation}

\begin{quotation}
2) The usage of the Wishart criteria leads to have a cardinality of features selected $ C^{W} $ according to the cardinality of the optimal clustering $ K^{*} $ computed with them, i.e., $ K^{*} \sim C^{W}$ is true. 
\end{quotation}

\begin{quotation}
3) The usage of Kaiser criteria does not lead to have similarity between the cardinality of features selected $ C^{K} $ with the cardinality of the optimal clustering $ K^{*} $ computed with them, i.e., $ K^{*} \sim C^{K} $ is false.
\end{quotation}

Therefore, the usage of Wishart distribution, and the derived Wishart limit $ \lambda_{max} $, helps to solve ill-posed clustering problem because it lows the matrix condition number and it satisfies the Ding and He theorem about the relation between true/optimal clusters and informative components, i.e.,

\begin{quotation}
4) PCA combined with Wishart criteria to select informative components works not only as data reduction method, but also as data clustering facilitation framework. 
\end{quotation}

Given the above conditions that suggest how to solve the clustering inverse problem with K-Means applied to data transformed via PCA, it is possible to affirm the following statement:

\begin{quotation}
5) Since there exists a similarity relation between number of informative components selected via Wishart criteria and the number of optimal clusters computed via K-Means clustering, i.e., $ K^{*} \sim C^{W}$., the cardinality of components selected as informative $ C^{W} $ suggests the cardinality of the true/optimal clusters to be found within the reduced dataset.
\end{quotation}


From the perspective of functional neuroimaging, there is an interesting result obtained with the clustering of task-related functional patterns within the 16 subjects, i.e.,

\begin{quotation}
6) There is no a direct relation between number of experimental stimuli and the number of temporal patterns found with data reduced via PCA and Wishart features selection. 
\end{quotation}

The aim of a temporal clustering it to find temporal patterns that are recurrent during the experimental session. The repetition of different stimuli is associated to four-to-seven temporal patterns, which are four-to-seven physiological states, i.e., independent brain reactions to different perceived stimuli, which are not the same of the experimental stimuli, but the representation and the elaboration of them by the brain once they have been presented. 
This \textit{distinguo} about perceived stimuli and experimental stimuli arises since it is plausible to assume that the temporal classes to be found within the fMRI dataset may be similar with the number of experimental stimuli. Since the stimuli are three for the task-blocks (famous faces, non-famous faces and scrambled faces) plus a rest-block, the true/optimal temporal patterns may be around three-to-four, whereas the clustering algorithm found four-to-seven classes\footnote{This statement is related to the K-Means applied to the data input reduced via PCA with components selected by Wishart criteria.}. This difference suggests that there is no a direct correspondence between the number of experimental stimuli and the number of different physiological status (i.e., the temporal patterns) acquired during the fMRI acquisition and detected by clustering algorithm.   

In conclusion, the relation between random matrix theory and clustering theory has salience in the statistical learning framework with application into the neuroscientific field. The proposed Wishart criteria to select informative components has been used also by authors as Chang et al \cite{chang2014cut} (MRI and fMRI studies) Burda \cite{burda2013collective} (fMRI study), Peyrache \cite{peyrache2010principal} (cell assemblies study). 

It is important to note that the limit defined by the Wishart distribution to select the informative eigenvectors has been used in a rigid manner with this study. The $ \lambda_{+} $ is the cut-off criteria for informative eigenvalues outside the distribution of the eigenvalues of gaussian random variables. However, from few years ago, in the random matrix theory domain, researcher found out also a distribution related to  $ \lambda_{+} $ (see works by Tracy and Widom \cite{tracy1994level,tracy1998correlation}, Nadal and Majumdar \cite{nadal2011simple} and the general research by Terence Tao on random matrix topics \cite{tao2012topics}).  Therefore, an extension of this work could be address the stochastic property of $ \lambda_{+} $ to draw a more fine association between clustering theory and random matrix theory.

Nevertheless, the results presented in this paper are related to the usage of K-Means clustering algorithms and its mathematical links with Principal Components Analysis (PCA). As mentioned in the introduction, clustering exists in the crisp and fuzzy form, and K-Means and PCA are methods belonged to the crisp category of unsupervised learning. So, a future work could be the investigation about the possible extension in the fuzzy domain of the theorem proposed by Ding and He for crisp clustering and crisp data reduction, using Fuzzy C-Means \cite{Bezdek1984,vergani2018soft} instead the K-Means \cite{lloyd1982least,Davies1979} and the Fuzzy PCA \cite{Luukka2011, Yang1999} instead the classic PCA \cite{Jolliffe2011}.

\begin{acknowledgement}
I would like to thank all the organizers of the Lake Como School of Advanced Studies that managed the School on Computational Methods For Inverse Problems in Imaging (May 21-25, 2018). I would also like to thank Prof. Elisabetta Binaghi and Prof. Marco Donatelli to support me in the acquisition of integrated knowledge about computational intelligence and numerical computing. 
\end{acknowledgement}

\bibliographystyle{unsrt}
\bibliography{bibliography}

\end{document}